\documentclass[12pt,leqno]{article}
\usepackage{amssymb,amscd,amsmath,amsthm}
\usepackage[hidelinks,hyperfootnotes=false]{hyperref}
\DeclareTextSymbolDefault{\textasteriskcentered}{OMS}
\hypersetup{
  pdftitle={The exceptional set of Goldbach numbers among polynomial values},
  pdfauthor={Ingo Alth\"ofer, Dijia Chen, Alberto Perelli and J\'anos Pintz},
  pdfsubject={Exceptional sets, corrected proofs and open problems}
}

\usepackage{xcolor}
\allowdisplaybreaks
\def\thefootnote{}

\def\beq{\begin{equation}}
\def\eeq{\end{equation}}

\theoremstyle{plain}
\newtheorem{Thm}{Theorem}

\theoremstyle{definition}
\newtheorem*{rema}{Remark}

\numberwithin{equation}{section}
\title{The exceptional set of Goldbach numbers among polynomial values}
\author{Ingo Alth\"ofer, Dijia Chen,\\
Alberto Perelli and J\'anos Pintz\thanks{Supported by the National Research
Development and Innovation Office, NKFIH, K~133819, K~147153 and
Excellence 151341.}}
\date{}

\begin{document}
\maketitle

\renewcommand{\thefootnote}{\arabic{footnote}}
\setcounter{footnote}{0}

\begin{abstract}
Let $\phi$ be a fixed polynomial with integer coefficients, positive leading
coefficient and degree $k\geq2$. We prove that the values $2\phi(n)$,
$1\leq n\leq N$, can be written as the sum of two primes, apart from at most
$N^{\alpha_k+o(1)}$ exceptional indices, where $\alpha_k<1$.
For quadratic polynomials the exponent is $1/2$, for cubic polynomials
$5/8$, for quartic polynomials $3/4$, and for $k>4$ it is
$1-1/(k+1)$. The proof combines a fourth-moment estimate over primes with
mean-value estimates over arbitrary subsets of polynomial sequences.
\end{abstract}

\section{Introduction}
\label{sec:1}

Let $\mathcal P$ denote the set of primes. Alth\"ofer and Pintz proved very
recently in a joint work \cite{PA2026} that the exceptional set in Goldbach's
problem,
\beq
\label{eq:1.1}
E(X):=\#\{n\in\mathbb N: n\leq X,\ 2\mid n,
\ n\neq p_1+p_2\text{ for all }p_1,p_2\in\mathcal P\},
\eeq
satisfies
\beq
\label{eq:1.2}
E(X)\ll_\varepsilon X^{7/12+\varepsilon}
\qquad(\varepsilon>0).
\eeq
This was improved a little later by H.~Li and Z.~Liu \cite{LL2026} to
\beq
\label{eq:1.3}
E(X)\ll_\varepsilon X^{1/2+\varepsilon}.
\eeq
The analytic arguments in these works use a fixed zero-free half-plane for
Dirichlet $L$-functions. The result of OpenAI quoted as Theorem~A below
supplies the required hypothesis. We use these external results in the
precise forms stated in Section~\ref{sec:2}.
Under the Generalized Riemann Hypothesis (GRH) for Dirichlet $L$-functions,
the classical result of Hardy and Littlewood \cite{HL1924} already gives
the exponent in \eqref{eq:1.3}.

Throughout the paper, $\phi\in\mathbb Z[x]$ is fixed, has degree $k\geq2$
and has positive leading coefficient. Put
\beq
\label{eq:1.4}
\begin{split}
\mathcal E_{k,\phi}(N)
&:=\{n\in\mathbb N: n\leq N,\ 2\phi(n)\neq p_1+p_2
\text{ for all }p_1,p_2\in\mathcal P\},\\
E_{k,\phi}(N)&:=\#\mathcal E_{k,\phi}(N).
\end{split}
\eeq
The finitely many indices for which $\phi(n)\leq0$ are included in this
definition and cause no difficulty. Since the values have size
$O_\phi(N^k)$, a direct application of \eqref{eq:1.3} gives no nontrivial
power bound for \eqref{eq:1.4} when $k\geq2$.

Perelli proved in 1996 \cite{Per1996} that, for every $A>0$,
\beq
\label{eq:1.5}
E_{k,\phi}(N)\ll_{\phi,A}N(\log N)^{-A}.
\eeq
This was improved by Br\"udern, Kawada and Wooley \cite{BKW2000} to
\beq
\label{eq:1.6}
E_{k,\phi}(N)\ll_\phi N^{1-c_0/k},
\eeq
with a small absolute constant $c_0>0$.
Our goal is to show a stronger explicit estimate for every $k\geq2$.

\begin{Thm}
\label{th:1}
Let $\phi\in\mathbb Z[x]$ be a fixed polynomial of degree $k\geq2$ with
positive leading coefficient. For every $\varepsilon>0$ and $N\geq2$,
\[
E_{k,\phi}(N)\ll_{\phi,\varepsilon}N^{\alpha_k+\varepsilon},
\]
where $\alpha_2=1/2$, $\alpha_3=5/8$, $\alpha_4=3/4$, and
$\alpha_k=1-1/(k+1)$ for $k>4$.
\end{Thm}

We use the real smooth weight $w$ specified in Section~2 of \cite{LL2026},
chosen in logarithmic coordinates from the Gevrey class of order~$2$.
In particular, $w\in C_c^\infty((0,1))$, $0\leq w\leq1$, and
$w(u)=1$ for $1/8\leq u\leq7/8$.
We use the circle method with
\beq
\label{eq:1.7}
S_w(\alpha):=\sum_{p\leq X}w(p/X)\log p\,e(p\alpha),
\qquad e(u)=e^{2\pi i u},\qquad L=\log X,
\eeq
and with polynomial sums of the form
\beq
\label{eq:1.8}
f_k(\alpha):=\sum_{1\leq n\leq N}e(2\phi(n)\alpha).
\eeq
Here $N$ measures indices and $X$ measures polynomial values. In the proof
we work with subsets for which $X/2<2\phi(n)\leq X$; on an index interval
$(N/2,N]$ all the required scales satisfy $X\asymp_\phi N^k$.

Write $\mathbb T=\mathbb R/\mathbb Z$ and define the major and minor arcs
on this circle by
\beq
\label{eq:1.9}
\begin{split}
\mathfrak M(P)
&=\bigcup_{q\leq P}\ \bigcup_{\substack{a\bmod q\\(a,q)=1}}
\left\{\alpha\in\mathbb T:
\left\|\alpha-\frac aq\right\|\leq\frac{P}{qX}\right\},\\
\mathfrak m(P)&=\mathbb T\setminus\mathfrak M(P).
\end{split}
\eeq
The norm denotes distance to the nearest integer. We choose $P=X^\kappa$
with a fixed $\kappa<1/2$, so that the individual major arcs are disjoint
for sufficiently large $X$. We usually omit $P$ from the notation.

For even $m\in[X/2,X]$, our goal is to obtain
\beq
\label{eq:1.10}
R(m):=\sum_{\substack{p+p'=m\\p,p'\leq X}}
w(p/X)w(p'/X)\log p\log p'>0
\eeq
when $m=2\phi(n)$, apart from a small exceptional set of indices.
Put
\beq
\label{eq:1.11}
\begin{split}
R_1(m)&=\int_{\mathfrak M}S_w(\alpha)^2e(-m\alpha)\,d\alpha,\\
R_2(m)&=\int_{\mathfrak m}S_w(\alpha)^2e(-m\alpha)\,d\alpha.
\end{split}
\eeq
Orthogonality gives
\beq
\label{eq:1.12}
R(m)=\int_0^1S_w(\alpha)^2e(-m\alpha)\,d\alpha=R_1(m)+R_2(m).
\eeq
For the simple weight $w=1$ on $[0,1]$, the classical goal is an
asymptotic formula of the type
\beq
\label{eq:1.13}
R_1^\star(m)=(1+o(1))\mathfrak S(m)m,
\eeq
where, for even positive $m$,
\beq
\label{eq:1.14}
\mathfrak S(m)=2C_0\prod_{\substack{p\mid m\\p>2}}
\left(1+\frac1{p-2}\right),\qquad
C_0=\prod_{p>2}\left(1-\frac1{(p-1)^2}\right)=0.66016\ldots.
\eeq
In the present weighted setting the corresponding goal is
\beq
\label{eq:1.15}
R_1(m)=(1+o(1))\mathfrak S(m)I_w(m),
\eeq
where
\beq
\label{eq:1.16}
I_w(m):=\sum_{r\in\mathbb Z}w(r/X)w((m-r)/X)
\geq m-X/4-1\gg X.
\eeq
Indeed, both weights are one for integers $r\in[X/8,m-X/8]$.
For our purposes it is sufficient to show, uniformly on the indicated
range of even $m$, that
\beq
\label{eq:1.17}
R_1(m)\geq c\mathfrak S(m)m
\eeq
with a fixed $c>0$.

\section{Auxiliary results for exponential sums over primes}
\label{sec:2}

Let $\theta$ denote the supremum of the real parts of the nontrivial zeros
of all primitive Dirichlet $L$-functions, including the Riemann zeta
function. We use the following result of OpenAI
\cite{OAI2026a,OAI2026b}.

\medskip\noindent
{\bf Theorem A.} {\it $\theta\leq7/8$.}

\medskip
For the weight fixed in Section~\ref{sec:1} and the arcs
\eqref{eq:1.9}, Li and Liu \cite{LL2026}, Theorem~1 and
Proposition~4.1, give the following two results.

\medskip\noindent
{\bf Theorem B.} {\it If $\theta<11/12$, there is a fixed
$\kappa=\kappa(\theta)\in(4/9,9/20)$ such that, with $P=X^\kappa$,}
\beq
\label{eq:2.1}
\int_{\mathfrak m}|S_w(\alpha)|^4\,d\alpha
\ll_{\theta,w}X^{5/2}(\log X)^{100}.
\eeq

\medskip\noindent
{\bf Theorem C.} {\it If $\theta<1$, then for every fixed
$\kappa\in(0,9/20)$ there is $c>0$ such that, with $P=X^\kappa$,
uniformly for even $m\in[X/2,X]$,}
\beq
\label{eq:2.2}
R_1(m)=\int_{\mathfrak M}S_w(\alpha)^2e(-m\alpha)\,d\alpha
=\mathfrak S(m)I_w(m)+O_{w,\theta,\kappa}(X^{1-c}).
\eeq

\begin{rema}
Pintz's explicit formula \cite{Pin2023} expresses the major-arc
contribution in terms of the pole at one and possible zeros near the line
$\operatorname{Re}s=1$. A fixed zero-free half-plane excludes the latter
zeros when $X$ is sufficiently large. For the precise range
$\kappa<9/20$ needed here we use Theorem~C. We choose the same $\kappa$
in Theorems~B and~C.
\end{rema}

Since $\mathfrak S(m)\geq2C_0>0$, Theorem~C and \eqref{eq:1.16}
imply \eqref{eq:1.17} for sufficiently large $X$. The integral $R_1(m)$
is real: the arcs are invariant under $\alpha\mapsto-\alpha$ on
$\mathbb T$, and $S_w(-\alpha)=\overline{S_w(\alpha)}$.
The proof below only requires the hypothesis $\theta<11/12$;
Theorem~A supplies it.

\section{Auxiliary results for exponential sums over polynomial values}
\label{sec:3}

To estimate the cardinality of the exceptional set, we use sums over its
subsets. For $\mathcal A\subseteq\{1,\ldots,M\}$ put
\beq
\label{eq:3.1}
K_{\mathcal A}(\alpha):=\sum_{m\in\mathcal A}e(2\phi(m)\alpha),
\qquad Z=\#\mathcal A.
\eeq
We require moment estimates that are uniform in the choice of
$\mathcal A$. The first follows from Wooley's weighted form of
Vinogradov's mean-value theorem \cite{Woo2019}, Corollary~1.4.

\medskip\noindent
{\bf Theorem D.} {\it Let $k\geq1$ and $s=k(k+1)/2$.
For every sequence $b_1,\ldots,b_M$ of complex numbers and every
$\varepsilon>0$,}
\beq
\label{eq:3.2}
\int_{\mathbb T^k}
\left|\sum_{m=1}^M b_m
e(\beta_1m+\beta_2m^2+\cdots+\beta_km^k)\right|^{2s}
d\boldsymbol\beta
\ll_{k,\varepsilon}M^\varepsilon
\left(\sum_{m=1}^M|b_m|^2\right)^s.
\eeq

\medskip
We use the following consequence.

\medskip\noindent
{\bf Corollary E.} {\it Let $\psi\in\mathbb Z[x]$ have degree
$k\geq1$ and nonzero leading coefficient. Put
\[
q=k(k+1),\qquad s=q/2,\qquad D=k(k-1)/2.
\]
For every $\mathcal A\subseteq\{1,\ldots,M\}$, $Z=\#\mathcal A$,
and every $\varepsilon>0$,}
\beq
\label{eq:3.3}
\int_{\mathbb T}
\left|\sum_{m\in\mathcal A}e(\alpha\psi(m))\right|^q\,d\alpha
\ll_{k,\varepsilon}M^{D+\varepsilon}Z^s.
\eeq
{\it The implied constant is independent of the coefficients of $\psi$.}

\begin{proof}
Let
\beq
\label{eq:3.FJ}
F(\boldsymbol\beta)
=\sum_{m\in\mathcal A}e(\beta_1m+\cdots+\beta_km^k),
\qquad J=\int_{\mathbb T^k}|F(\boldsymbol\beta)|^{2s}
\,d\boldsymbol\beta.
\eeq
By Theorem~D, $J\ll_{k,\varepsilon}M^\varepsilon Z^s$.
For $\mathbf h=(h_1,\ldots,h_k)\in\mathbb Z^k$, let
$N(\mathbf h)$ count the ordered pairs
$(\mathbf x,\mathbf y)\in\mathcal A^s\times\mathcal A^s$ such that
\[
\sum_{i=1}^s(x_i^j-y_i^j)=h_j\qquad(1\leq j\leq k).
\]
Orthogonality gives
\[
N(\mathbf h)=\int_{\mathbb T^k}|F(\boldsymbol\beta)|^{2s}
e(-\boldsymbol\beta\cdot\mathbf h)\,d\boldsymbol\beta,
\qquad 0\leq N(\mathbf h)\leq J.
\]
Write $\psi(m)=a_km^k+\cdots+a_0$. Since there are equally many
$x_i$ and $y_i$, the constant coefficient cancels. A second application
of orthogonality yields
\beq
\label{eq:3.4}
\int_{\mathbb T}
\left|\sum_{m\in\mathcal A}e(\alpha\psi(m))\right|^{2s}\,d\alpha
=\sum_{\substack{\mathbf h\in\mathbb Z^k\\a_1h_1+\cdots+a_kh_k=0}}
N(\mathbf h).
\eeq
A nonzero summand has $|h_j|\leq sM^j$. Once
$h_1,\ldots,h_{k-1}$ are fixed, the relation in \eqref{eq:3.4}
determines at most one integer $h_k$, because $a_k\neq0$.
Thus there are at most
\[
\prod_{j=1}^{k-1}(2sM^j+1)\ll_k M^{1+\cdots+(k-1)}=M^D
\]
potentially nonzero summands. For $k=1$ the product is empty and equals
one. Bounding each summand by $J$ proves \eqref{eq:3.3}.
\end{proof}

For $k=3$ and $k=4$ we use the following result of Hughes and Wooley
\cite{HW2022}, Theorem~4.1, equation~(4.2). It is valid for every
degree at least three.

\medskip\noindent
{\bf Theorem F.} {\it If $\psi\in\mathbb Z[x]$ has degree at least
three, $\mathcal A\subseteq\{1,\ldots,M\}$ and $Z=\#\mathcal A$, then}
\beq
\label{eq:3.5}
\int_0^1\left|\sum_{m\in\mathcal A}e(\alpha\psi(m))\right|^8\,d\alpha
\ll_{\psi,\varepsilon}M^{1+\varepsilon}Z^4.
\eeq

\medskip
For degree two an elementary divisor argument suffices.

\medskip\noindent
{\bf Lemma G.} {\it If $\psi\in\mathbb Z[x]$ has degree two,
$\mathcal A\subseteq\{1,\ldots,M\}$ and $Z=\#\mathcal A$, then}
\beq
\label{eq:3.6}
\int_0^1\left|\sum_{m\in\mathcal A}e(\alpha\psi(m))\right|^4\,d\alpha
\ll_{\psi,\varepsilon}M^\varepsilon Z^2.
\eeq

\begin{proof}
Write $\psi(m)=am^2+bm+c$, with $a\neq0$, and put
\[
r(t)=\#\{(x,y)\in\mathcal A^2:\psi(x)-\psi(y)=t\}.
\]
Orthogonality gives
\beq
\label{eq:3.7}
\int_0^1\left|\sum_{m\in\mathcal A}e(\alpha\psi(m))\right|^4\,d\alpha
=\sum_{t\in\mathbb Z}r(t)^2,
\qquad \sum_t r(t)=Z^2.
\eeq
For $t\neq0$ the factorization
\beq
\label{eq:3.8}
(x-y)\bigl(a(x+y)+b\bigr)=t
\eeq
gives at most $2d(|t|)$ choices for the signed divisor $z=x-y$.
Once $z$ is fixed, the other factor is fixed, and
\beq
\label{eq:3.9}
a(2x-z)+b=t/z
\eeq
determines at most one $x$. Since $|t|\ll_\psi M^2$, the divisor
bound, with a suitably smaller exponent, gives
$r(t)\ll_{\psi,\varepsilon}M^\varepsilon$ for $t\neq0$.
For $t=0$, either $x=y$ or $a(x+y)+b=0$, whence $r(0)\leq2Z$.
Consequently,
\[
\sum_t r(t)^2
\leq4Z^2+O_{\psi,\varepsilon}(M^\varepsilon)
\sum_{t\neq0}r(t)
\ll_{\psi,\varepsilon}M^\varepsilon Z^2.
\]
\end{proof}

\begin{rema}
The diagonal solutions give a lower bound $Z^2$ for the moment in
\eqref{eq:3.6}. Thus its dependence on $Z$ is sharp, apart from the
factor $M^\varepsilon$.
\end{rema}

\section{Proof of Theorem 1}
\label{sec:4}

Fix $\phi$ and let $N$ be sufficiently large in terms of its coefficients.
For indices $N/2<n\leq N$, positivity of the leading coefficient gives
\[
c_1N^k\leq2\phi(n)\leq c_2N^k
\]
with constants $c_1,c_2>0$ depending only on $\phi$. This value interval
is covered by $O_\phi(1)$ consecutive dyadic intervals $(X/2,X]$,
all with $X\asymp_\phi N^k$. For one such interval put
\[
\mathcal A=
\{n\in\mathcal E_{k,\phi}(N):N/2<n\leq N,\ X/2<2\phi(n)\leq X\},
\qquad Z=\#\mathcal A.
\]
At this scale choose $P=X^\kappa$ with the fixed $\kappa$ of Theorem~B.
Theorem~C applies to the same arcs. For every $n\in\mathcal A$,
$R(2\phi(n))=0$, and therefore
\beq
\label{eq:4.1}
R_2(2\phi(n))=-R_1(2\phi(n))
=-(1+o(1))\mathfrak S(2\phi(n))I_w(2\phi(n))\leq-cX
\eeq
for a fixed $c>0$, uniformly in $n$.
Define
\beq
\label{eq:4.2}
K(\alpha)=\sum_{n\in\mathcal A}e(-2\phi(n)\alpha).
\eeq
Summing \eqref{eq:4.1}, interchanging the finite sum and the integral,
and taking absolute values gives
\beq
\label{eq:4.3}
ZX\ll\left|\int_{\mathfrak m}S_w(\alpha)^2K(\alpha)\,d\alpha\right|
\leq\int_{\mathfrak m}|S_w(\alpha)|^2|K(\alpha)|\,d\alpha.
\eeq

Suppose that, for some $q\geq4$ and $b\geq0$, we have
\beq
\label{eq:4.4}
\left(\int_0^1|K(\alpha)|^q\,d\alpha\right)^{1/q}
\ll_{\phi,\varepsilon}N^{b/q+\varepsilon}Z^{1/2}.
\eeq
All results of Section~\ref{sec:3} apply with $\psi=-2\phi$.
Lemma~G gives $q=4$, $b=0$ for $k=2$.
Theorem~F gives $q=8$, $b=1$ for $k=3,4$.
Corollary~E gives $q=k(k+1)$, $b=k(k-1)/2=q/2-k$ in general.

We write
\beq
\label{eq:4.5}
|S_w(\alpha)|^2
=\bigl(|S_w(\alpha)|^4\bigr)^{1/q}
\bigl(|S_w(\alpha)|^2\bigr)^{1-2/q}
\eeq
and apply H\"older's inequality to these two factors and $|K|$, with
exponents $q$, $q/(q-2)$ and $q$, respectively.
Parseval's identity and the elementary estimate
$\sum_{p\leq X}\log p\ll X$ give
\[
\int_0^1|S_w(\alpha)|^2\,d\alpha
=\sum_{p\leq X}w(p/X)^2(\log p)^2\ll XL.
\]
Together with Theorem~B and \eqref{eq:4.4}, this yields
\begin{align}
\label{eq:4.6}
ZX
&\ll
\left(\int_{\mathfrak m}|S_w|^4\right)^{1/q}
\left(\int_0^1|S_w|^2\right)^{1-2/q}
\left(\int_0^1|K|^q\right)^{1/q}\\
&\ll_{\phi,\varepsilon}
X^{1+1/(2q)}L^{1+98/q}N^{b/q+\varepsilon}Z^{1/2}.
\nonumber
\end{align}
If $Z=0$ there is nothing to prove. Otherwise divide by $XZ^{1/2}$
and square. Since $X\asymp_\phi N^k$, and since $\varepsilon$ may be
chosen arbitrarily small, we obtain
\beq
\label{eq:4.7}
Z\ll_{\phi,\varepsilon}N^{(k+2b)/q+\varepsilon}.
\eeq
Thus the quantity to minimize among the available moment estimates is
$(k+2b)/q$.

For $k=2$, the choices $q=4$, $b=0$ give $\alpha_2=1/2$.
For $k=3,4$, the choices $q=8$, $b=1$ give
$\alpha_k=(k+2)/8$.
For $k>4$, the choices $q=k(k+1)$, $2b=q-2k$ give
$\alpha_k=1-k/q=1-1/(k+1)$.
Summing over the $O_\phi(1)$ value intervals proves the required bound
for the exceptional indices in $(N/2,N]$. Summing next over the index
intervals $(N/2^{j+1},N/2^j]$, and absorbing the finitely many initial
indices into the implied constant, proves Theorem~\ref{th:1}.
\hfill$\square$

\section{Application to other thin sequences}
\label{sec:5}

Theorem~\ref{th:1} gives related thin sequences containing only finitely
many Goldbach exceptional numbers. For example, consider the values
$2m^k$ with $m=\lfloor a(n+\beta)^r\rfloor$, where $a>0$ is fixed.
For almost every $\beta\in[0,1]$ there are only finitely many exceptions
when $r>1/(1-\alpha_k)$. Thus it suffices to take $r>2$ for $k=2$,
$r>8/3$ for $k=3$, $r>4$ for $k=4$, and $r>k+1$ for $k>4$.
No parity restriction on $m$ is needed.

Put
\[
\mathcal B_k=\{m\in\mathbb N:2m^k\neq p+q
\text{ for all primes }p,q\}.
\]
Assume, for a fixed $0\leq c_k<1$, that for every $\varepsilon>0$
\beq
\label{eq:5.1}
\#(\mathcal B_k\cap[1,Y])\ll_{k,\varepsilon}Y^{c_k+\varepsilon}
\qquad(Y\geq2).
\eeq

\medskip\noindent
{\bf Proposition.} {\it Fix $a>0$ and $r>1/(1-c_k)$, with
$0\leq c_k<1$ as in \eqref{eq:5.1}. For almost every
$\beta\in[0,1]$, every sufficiently large $n$ satisfies}
\beq
\label{eq:5.2}
2\lfloor a(n+\beta)^r\rfloor^k=p_n+q_n
\qquad(p_n,q_n\text{ prime}).
\eeq
{\it Moreover, for
$\mathcal A_{a,r,\beta}=\{2\lfloor a(n+\beta)^r\rfloor^k:n\geq1\}$,
the counting asymptotic}
\[
\#(\mathcal A_{a,r,\beta}\cap[1,X])
\sim\left(\frac{X}{2a^k}\right)^{1/(kr)}
\]
{\it holds for every fixed $\beta\in[0,1]$.}

\begin{proof}
Put $m_n(\beta)=\lfloor a(n+\beta)^r\rfloor$ and
$B_{k,n}=\{\beta\in[0,1]:m_n(\beta)\in\mathcal B_k\}$.
The assumption on $r$ implies $r>1$.
For $T\leq n<2T$, the derivative of $a(n+\beta)^r$ with respect to
$\beta$ is $\asymp_{a,r}T^{r-1}$. Hence, for every integer $b$,
\beq
\label{eq:5.3}
\operatorname{meas}\{\beta\in[0,1]:m_n(\beta)=b\}
\ll_{a,r}T^{1-r}.
\eeq
The intervals $[an^r,a(n+1)^r]$ have disjoint interiors. Their lengths
tend to infinity since $r>1$. For sufficiently large $T$, each unit
interval $[b,b+1)$ therefore meets at most two of them, apart from
irrelevant endpoints. It follows from \eqref{eq:5.1} and
\eqref{eq:5.3} that
\begin{align}
\label{eq:5.4}
\sum_{T\leq n<2T}\operatorname{meas}(B_{k,n})
&\ll_{a,r}T^{1-r}
\#\bigl(\mathcal B_k\cap[1,C_{a,r}T^r]\bigr)\\
&\ll_{a,r,k,\varepsilon}T^{1-r+ r(c_k+\varepsilon)}
=T^{1-r(1-c_k-\varepsilon)}.
\nonumber
\end{align}
Choose $\varepsilon>0$ so that $r(1-c_k-\varepsilon)>1$.
Summation over $T=2^j$ gives
$\sum_n\operatorname{meas}(B_{k,n})<\infty$.
The first Borel--Cantelli lemma proves \eqref{eq:5.2}; no independence
assumption is required. Also $m_n(\beta)\geq1$ for all sufficiently
large $n$, so that the omission of $0$ from $\mathcal B_k$ is harmless.
Finally, $m_n(\beta)\sim an^r$ and the values are eventually strictly
increasing. Inverting $2m_n(\beta)^k\sim2a^k n^{kr}$ proves the counting
asymptotic.
\end{proof}

\begin{rema}
The measure argument does not provide an explicit admissible $\beta$.
For example, can one exhibit a particular $\beta\in[0,1]$ for which
\beq
\label{eq:5.5}
2\lfloor(n+\beta)^3\rfloor^2
\eeq
is the sum of two primes for every sufficiently large $n$?
Section~\ref{sec:5} proves this for almost every $\beta$, by taking
$k=2$, $a=1$ and $r=3$. The particularly natural choice $\beta=0$
would assert that every sufficiently large number $2n^6$ is Goldbach,
a substantially stronger conclusion than the exceptional-set estimates
proved in this paper. We do not know any triple $(a,k, r)$ with $a>0$, 
integer $k>0, r>0$ for which we could specify a value
$\beta\in (0,1]$ such that the expression \eqref{eq:5.2}  would be a sum of two
primes for every value of n.
\end{rema}

\medskip\noindent
{\bf AI assistance.}
OpenAI Astra was used for careful auditing. Its checks concerned the
proof, references and notation. It also assisted with the revision and
with formulating the above remark.

\smallskip\noindent
{\small
Ingo Alth\"ofer\\
Faculty of Mathematics and Computer Science\\
Friedrich-Schiller University Jena, D-07740 Jena, Germany\\
E-mail: ingo.althoefer@uni-jena.de

\smallskip\noindent
Dijia Chen\\
Department of Mathematics, University of Wisconsin--Madison\\
Madison, WI 53706, USA\\
E-mail: dchen426@wisc.edu

\smallskip\noindent
Alberto Perelli\\
Dipartimento di Matematica, Universit\`a di Genova\\
Via Dodecaneso 35, 16146 Genova, Italy\\
E-mail: alberto.perelli@unige.it

\smallskip\noindent
J\'anos Pintz\\
HUN-REN Alfr\'ed R\'enyi Institute of Mathematics\\
Re\'altanoda u.~13--15, H-1053 Budapest, Hungary\\
E-mail: pintz@renyi.hu}

\end{document}